\documentclass[reqno,12pt]{amsart}
\usepackage{bbm}
\usepackage{amsmath}
\usepackage{txfonts}
\usepackage{stmaryrd}
\usepackage{amssymb}
\usepackage{amsfonts}
\usepackage{times}
\usepackage{mathrsfs}
\usepackage{galois}

\usepackage{amsthm}
\usepackage{enumerate}

\usepackage{framed}
\usepackage{lipsum}
\usepackage{color}

\newcommand{\DD}{{\mathbb D}}
\def\B{\mathcal{B}}
\def\R{\mathcal{R}}

\def\D{\mathbb{D}}
\def\N{\mathbb N}

\def\vp{\varphi}

\def\p{{\prime}}

\def\msk{\medskip}
\def\ol{\overline}

\def\BB{ \mathbb{B}}
\def\CC{ \mathbb{C}}
\def\mu{\mathcal \beta}

\begin{document}
\title[ Weighted Composition Operators  ]{ Weighted Composition Operators  from $H^\infty$ to the Bloch  Space in the Unit Ball of $\CC^n$}

\author{ Juntao Du and  Songxiao Li$^\dagger$  }

\address{Juntao Du\\ Faculty of Information Technology, Macau University of Science and Technology, Avenida Wai Long, Taipa, Macau.}
\email{ jtdu007@163.com  }
\address{Songxiao Li\\ Institute of Fundamental and Frontier Sciences, University of Electronic Science and Technology of China,
610054, Chengdu, Sichuan, P.R. China\newline
Institute of Systems Engineering, Macau University of Science and Technology, Avenida Wai Long, Taipa, Macau. } \email{jyulsx@163.com}

\begin{abstract}  The boundedness and compactness of weighted composition operators from $H^\infty$ to the Bloch  space  in the unit ball of $\CC^n$ are investigated in this paper. In particular,   some new characterizations for the boundedness and  the essential norm of weighted composition operators are given.
\thanks{$\dagger$ Corresponding author.}
\thanks{This  project was partially supported by NSF of China (No.11471143 and No. 11720101003). }
\end{abstract}

\maketitle

\section{Introduction}

Let $\BB$ be the open unit ball of $\CC^n$ and $\partial \BB$   the boundary of $\BB$. When $n=1$, $\BB$ is the open unit disk $\D$ in the complex plane. Let $H(\BB)$ denote the space of all holomorphic functions on $\BB$. For $f\in H(\BB)$, the  radial derivative and complex gradient of  $f$ at $z$ will be denoted by $\R f(z)$ and  $\nabla f(z)$, respectively. That is,
$$\R f(z)=\sum_{j=1}^n z_j\frac{\partial f}{\partial z_j}(z),\mbox{  and }
\,\,\nabla f(z)=\Big(\frac{\partial f}{\partial z_1}(z),\frac{\partial f}{\partial z_2}(z),\cdots,
\frac{\partial f}{\partial z_n}(z)\Big).$$

An $f\in H(\BB)$ is said to belong to the Bloch space, denoted by
 $\B=\B(\BB)$, if
$$ \|f\|_{\beta}:=\sup\limits_{z\in \BB} (1-|z|^2)|\R f(z)|<\infty.$$
The space $\B$ is a Banach space with the norm $\|f\|_{\B}=|f(0)|+\|f\|_\mu $.
From \cite{Zkh2005}, we see that $\|f\|_\beta \approx \sup_{z\in\BB}(1-|z|^2)|\nabla f(z)|.$
In  \cite{Tr1980blms},  Timoney  proved that
\begin{equation}\label{1222-3}
\|f\|_{\beta}\approx \sup\left\{ \frac{|\langle \nabla f(z),\ol{v} \rangle|}{H_z(v,v)^\frac{1}{2}} :z\in\BB,v\in\CC^n\backslash\{0\}  \right\}.
\end{equation}
Here $H_z(v,v)$ is the Bergman metric defined by
$$H_z(v,v)=\frac{n+1}{2}\frac{(1-|z|^2)|v|^2+|\langle v,z\rangle|^2}{(1-|z|^2)^2}, ~~z\in\BB,v\in\CC^n\backslash\{0\}  .$$
See \cite{Zkh2005} for
more information of the Bloch space $\B$ on the unit ball.

We   use $H^\infty$ to denote the  space of bounded holomorphic functions in $\BB$. That is, $f\in H^\infty$ if and only if $f\in H(\BB)$ and
$$\|f\|_{\infty}:=\sup_{z\in\BB} |f(z)|<\infty.$$
It is well known that $H^\infty$ is a Banach space and a subset of $\B$. Moreover (see \cite{Zkh2005}),
\begin{equation}\label{1222-1}
\|f\|_\B \leq  \|f\|_{\infty}. \end{equation}

Let $u\in H(\BB)$ and $\vp(z)=(\vp_1(z),\vp_2(z),\cdots,\vp_n(z))$ be a holomorphic self-map of $\BB$. The weighted composition operator, denoted by  $uC_\vp$,  is defined by
$$(uC_\vp f)(z)=u(z)f(\vp(z)),\,\,f\in H(\BB).$$
When $u=1$, $uC_\vp$ is the composition operator, denoted by $C_\vp$. It is important to give function theoretic description of when $u$ and $\vp$ induce a bounded or compact weighted composition operator on various function spaces (see \cite{CccMbd1995}).

 In the setting of the unit disk, it is   well known that the operator  $C_\vp$ is bounded on $\B$ for any analytic self-map $\vp$ in $\D$ by  Schwarz-Pick Lemma. The compactness of the composition operator on $\B$ was characterized in \cite{MkMa1995tams}.  In \cite{WlZdcZkh2009pams}, Wulan, Zheng and Zhu  showed that $C_\varphi$ is compact on $\B$ if and only if $
\lim\limits_{k\to\infty} \|\varphi^k \|_\B=0.$ This method has been used to describe the boundedness and compactness of $uC_\varphi$ on some function spaces, see \cite{Cf2013cejm,Djn2012jmaa,HqhLsxWh2017cvee,LxsLsx2017ieot} for example.
In \cite{Os2001tjm}, Ohno  characterized the boundedness and compactness of the operator  $uC_\varphi:H^\infty\to\B$. Colonna, motivated by \cite{WlZdcZkh2009pams}, gave another characterization for the boundedness and compactness of the operator $uC_\varphi:H^\infty\to\B$ in \cite{Cf2013cejm}.  Hu, Li and Wulan, based on the work of  Ohno and  Colonna,  gave some estimates for the essential norm  of the operator  $uC_\varphi:H^\infty\to\B$ in \cite{HqhLsxWh2017cvee}.  Moreover, they gave a new characterization for the boundedness and compactness of $uC_\varphi:H^\infty\to\B$ in \cite{HqhLsxWh2017cvee}.

 In the setting of the unit ball, Shi and  Luo studied   composition operators on the Bloch space in \cite{SjhLl2000ams}. In \cite{Djn2012jmaa}, Dai gave several new characterizations for the compactness of the composition operator  on the Bloch space, which extended  the main result in \cite{WlZdcZkh2009pams} to the unit ball. Li and Stevi\'c characterized the boundedness and compactness of $uC_\varphi:H^\infty\to\B$ in \cite{LsxSs2008tjm} (see \cite{LsxSs2007aaa} for the setting of polydisk). Zhang and Chen gave two characterizations of the boundedness and compactness of $uC_\varphi:H^\infty\to\B$ in \cite{ZmzChh2009ams}.  For example, they showed that $uC_\vp:H^\infty\to\B$ is bounded
 if and only if $u\in\B$  and
  $$ \sup_{z\in\BB} (1-|z|^2)|u(z)|\sqrt{H_{\vp(z)}(\vp^\p(z) z,\vp^\p(z) z)}<\infty.$$

In this paper, motivated by \cite{HqhLsxWh2017cvee} and \cite{ZmzChh2009ams}, we investigate the boundedness, compactness and essential norm of $uC_\vp:H^\infty\to\B$ on the unit ball. That is, we will give some new characterizations for the boundedness, compactness and essential norm of $uC_\vp:H^\infty\to\B$. These extend the results in \cite{HqhLsxWh2017cvee} to the unit ball. Moreover, the method we used here is completely different from
\cite{HqhLsxWh2017cvee}.

Recall that the essential norm of a bounded linear operator $T:X\to Y$,  denoted by $\|T\|_{e,X\to Y}$, is defined
as the distance from $T$ to the space of compact operators  from $X$ to $Y$. That is,
$$\|T\|_{e,X\to Y}=\inf\{\|T-K\|_{X\to Y}: \, \, ~~~K \mbox{ is a compact operator from } X \mbox{ to } Y\}.$$

Constants are denoted by $C$, they are positive and
may differ from one occurrence to the next. We say that $A\lesssim B$ if there exists a constant $C$ such that $A \leq CB$. The symbol $A\approx  B$ means that $A\lesssim B\lesssim A$.

 \section{Boundedness of $uC_\vp:H^\infty\to\B  $}
Before we state the main result and the proof in this section, we state some notations and preliminary results. Let $\vp^\p(z)$ be the  Jacobian matrix of $\vp$, that is
 $$\vp^\p(z)=\left(
 \begin{array}{cccc}
\frac{\partial \vp_1}{\partial z_1} &\frac{\partial \vp_1}{\partial z_2} &\cdots&\frac{\partial \vp_1}{\partial z_n}\\
\frac{\partial \vp_2}{\partial z_1} &\frac{\partial \vp_2}{\partial z_2} &\cdots&\frac{\partial \vp_2}{\partial z_n}\\
\cdots                              &\cdots                               &\cdots     &\cdots\\
\frac{\partial \vp_n}{\partial z_1} &\frac{\partial \vp_n}{\partial z_2} &\cdots&\frac{\partial \vp_n}{\partial z_n}
 \end{array}
 \right).$$
 Therefore
 $$
 \nabla (f(\vp(z)))=(\nabla f)(\vp(z))\vp^\p(z),\mbox{ and   }\,\, \R (f(\vp(z)))=(\nabla f)(\vp(z))\vp^\p(z)z.
 $$
 Here and henceforth,  we do not distinguish the row vector and column vector, that is, we always  admit the vectors   have the proper forms in the expressions.

For $a\in \BB\backslash\{0\}$,  the  automorphism of $\BB$ is defined by
$$\phi_a(z)=\frac{a-P_a z-s_aQ_a z}{1-\langle z,a\rangle}, z\in\BB,$$
where $s_a=\sqrt{1-|a|^2}$,
$$P_a z=\frac{\langle z,a\rangle}{|a|^2}a,\,\, Q_a z=z-\frac{\langle z,a\rangle}{|a|^2}a,\,z\in\BB.$$
When $a=0$, set $\phi_a(z)=-z$.
Let $\phi_a(z)=(\phi_{a,1}(z),\phi_{a,2}(z),\cdots,\phi_{a,n}(z))$. We have
$$\{\phi_{a,i}(z)\}_{i=1}^n\subset H^\infty,\mbox{ and } \sum_{i=1}^n |\phi_{a,i}(z)|^2<1.$$

\noindent{\bf Lemma 1.} \cite[Lemma 2.1]{Djn2012jmaa} {\it Let $a\in \BB\backslash\{0\}$. Then
$$|\phi_a(z)-a|=\frac{\sqrt{(1-|a|^2)(|z|^2-|\langle z,a \rangle|^2)}}{|1-\langle z,a \rangle|}
 $$
and
$$|\phi_a^\p(a)z|=\frac{\sqrt{(1-|a|^2)|z|^2+ |\langle z,a \rangle|^2}}{1-|a|^2}.$$
}

\noindent{\bf Lemma 2.} \cite{li} {\it Let $u,f\in H(\BB)$. Then
$$\R I_u(f)(z)=u(z)\R f(z),\,\,\R J_u(f)(z)= \R u(z)f(z), ~~~\,~\, z\in \BB.$$
Here
$$I_u(f)(z)=\int_0^1 u(tz)(\R f)(tz)\frac{dt}{t},\,\,J_u(f)(z)=\int_0^1 (\R u)(tz)f(tz)\frac{dt}{t}.$$
}

\noindent{\bf Theorem 1.} {\it Suppose $u\in H(\BB)$ and $\vp$ is a holomorphic self-map of $\BB$.
Then the following statements are equivalent.
\begin{enumerate}[(i)]
  \item $uC_\vp:H^\infty\to\B$ is bounded.
  \item $u\in\B$ and $$M_1:=\sup_{k\in\N}\sup_{\xi\in\partial \BB}  \|u\langle \varphi,\xi\rangle^k\|_{\beta} <\infty.$$
  \item $u\in\B$ and
  $$M_2:=\sup_{z\in\BB} (1-|z|^2)|u(z)|\sqrt{H_{\vp(z)}(\vp^\p(z) z,\vp^\p(z) z)}<\infty.$$
  \item $u\in\B$ and
  $$M_3:=\sup_{1\leq i\leq n}\sup\limits_{w\in\BB}\|uC_\vp\phi_{\vp(w),i}\|_{\beta}<\infty.$$
  \item $u\in\B$ and
  $$M_4:=\sup_{k\in\N}\sup_{\xi\in\partial \BB} \|I_u(\langle \varphi,\xi\rangle^k)\|_{\beta}<\infty.$$
\end{enumerate}
}
\begin{proof} Obviously,  $uC_\vp:H^\infty\to\B$ is bounded if and only if
$$\|uC_\vp f\|_{\mu}\lesssim \|f\|_{\infty}, ~~~\forall f\in H^\infty.$$

 {\bf ({\it i})$\Rightarrow$ ({\it ii})}.  This implication  is obvious since $\|f_{k,\xi}\|_{\infty}=1$.  Here and henceforth, $f_{k,\xi}(z)=\langle z,\xi\rangle^k$, $z\in \BB$, $\xi\in\partial \BB$.

{\bf ({\it ii})$\Rightarrow$ ({\it iii})}. When $k\geq 1$, since
$$\R(f_{k,\xi}\comp\vp)(z)=k \langle \vp(z),\xi\rangle^{k-1} \langle \vp^\p(z) z,\xi  \rangle$$
and
\begin{equation*}
|\R (uC_\vp f_{k,\xi})(z)|
\geq  |u(z)||\R(f_{k,\xi}\comp\vp)(z)|-|\R u(z)||f_{k,\xi}(\vp(z))|  ,
\end{equation*}
we obtain
\begin{equation}\label{1216-1}
\sup_{z\in\BB}\sup_{k\in\N}\sup_{\xi\in\partial \BB}k(1-|z|^2)|u(z)|\left| \langle \vp(z),\xi\rangle^{k-1}  \langle \vp^\p(z) z,\xi  \rangle \right| \leq 2M_1.
\end{equation}
After a calculation, we have
{\small
\begin{equation}\label{1216-2}
\sqrt{H_{\vp(z)}(\vp^\p(z) z,\vp^\p(z) z) }
\lesssim \left( \frac{(1-|\vp(z)|^2)^\frac{1}{2}|\vp^\p(z)z|}{1-|\vp(z)|^2}+\frac{|\langle  \vp^\p(z)z,\vp(z) \rangle|}{1-|\vp(z)|^2}\right).
\end{equation}
}

Let $\xi_i$ be the vector in which the  $i$-th component is 1 and the others are 0, $i=1,2,\cdots,n$.
By letting $k=1$, from (\ref{1216-1}), we have
\begin{equation}\label{1220-2}
\sup_{z\in\BB}(1-|z|^2)|u(z)||\vp^\p(z) z|
=\sup_{z\in\BB}(1-|z|^2)|u(z)|\left(\sum_{i=1}^n |\langle \vp^\p(z) z,\xi_i \rangle|^2\right)^\frac{1}{2} \lesssim M_1.
\end{equation}

When $|\vp(z)|\leq \frac{1}{2}$, from (\ref{1216-2}) we have
\begin{equation*}
 (1-|z|^2)|u(z)|\sqrt{H_{\vp(z)}(\vp^\p(z) z,\vp^\p(z) z)} \lesssim  (1-|z|^2)|u(z)||\vp^\p(z) z|.
\end{equation*}
Therefore,
\begin{equation}\label{1216-3}
 \sup_{|\vp(z)|\leq \frac{1}{2}}(1-|z|^2)|u(z)|\sqrt{H_{\vp(z)}(\vp^\p(z) z,\vp^\p(z) z)}
 \lesssim   M_1.
\end{equation}

Next, we will prove that
\begin{equation}\label{0101-1}
 \sup_{|\vp(z)|> \frac{1}{2}}(1-|z|^2)|u(z)|\sqrt{H_{\vp(z)}(\vp^\p(z) z,\vp^\p(z) z)}\lesssim M_1.
\end{equation}
Assume that there exists $k\in\N$ such that $k\geq 3 $ and $1-\frac{1}{k-1}\leq |\vp(z)| < 1-\frac{1}{k}$.
It is easy to see that 
\begin{equation}\label{0102-1}
     |\vp(z)|^{2(k-1)}\approx 1  ,  \mbox{ and } k(  1-|\vp(z)|^2)\approx 1.
\end{equation}
Therefore, by letting  $\tau=\vp(z)$, we have
\begin{equation}\label{1219-6}
\frac{|\langle  \vp^\p(z)z,\vp(z) \rangle|}{1-|\vp(z)|^2}
\approx  k {|\langle \vp(z),\tau \rangle|^{k-1}|\langle  {\vp^\p(z)z,\tau \rangle|}}.
\end{equation}
By (\ref{1216-1}), we get
\begin{eqnarray}  \sup_{1-\frac{1}{k-1}\leq |\vp(z)| < 1-\frac{1}{k} } \frac{(1-|z|^2)|u(z)||\langle  \vp^\p(z)z,\vp(z) \rangle|}{1-|\vp(z)|^2} \lesssim   M_1  \nonumber
\end{eqnarray}
and hence
\begin{eqnarray}\label{1216-4}
& & \sup_{|\vp(z)|> \frac{1}{2}} \frac{(1-|z|^2)|u(z)||\langle  \vp^\p(z)z,\vp(z) \rangle|}{1-|\vp(z)|^2} \nonumber\\
&\lesssim &  \sup_{k \geq 3} \sup_{1-\frac{1}{k-1}\leq |\vp(z)| < 1-\frac{1}{k} } \frac{(1-|z|^2)|u(z)||\langle  \vp^\p(z)z,\vp(z) \rangle|}{1-|\vp(z)|^2} \lesssim   M_1.
\end{eqnarray}
By projection theorem, there exists a $\eta(z)\in \partial \BB$ such that $\langle \vp(z),\eta(z)\rangle =0. $ Then
\begin{equation}\label{1222-2}
\vp^\p(z)z=p\vp(z)+q \eta(z),
\end{equation}
where
\begin{equation*}
p=\frac{\langle \vp^\p(z) z,\vp(z)\rangle}{|\vp(z)|^2}, q=\langle \vp^\p(z) z,\eta(z)\rangle.
\end{equation*}
Since $|\vp(z)|>\frac{1}{2}$, we get
\begin{equation}\label{1217-1}
|\vp^\p(z)z|^2 \leq 4|\langle \vp^\p(z) z,\vp(z)\rangle|^2+|\langle \vp^\p(z) z,\eta(z)\rangle|^2.
\end{equation}
Let $\zeta(z)=\vp(z)+\sqrt{1-|\vp(z)|^2}\eta(z)$. By (\ref{1222-2}), we have
\begin{equation}\label{1222-4}
|\zeta(z)|=1,\,\, \langle \vp(z),\zeta(z)\rangle =|\vp(z)|^2,
\end{equation}
and
$$\langle \vp^\p(z)z,\zeta(z)\rangle =\langle  \vp^\p(z)z,\vp(z) \rangle +\sqrt{1-|\vp(z)|^2}\langle \vp^\p(z) z,\eta(z)\rangle.$$
Then,
\begin{eqnarray}
\sqrt{1-|\vp(z)|^2}|\langle \vp^\p(z) z,\eta(z)\rangle|
&\leq&   |\langle \vp^\p(z) z,\vp(z)\rangle| +|\langle \vp^\p(z) z,\zeta(z)\rangle| .   \label{1217-2}
\end{eqnarray}
From (\ref{1217-1}), (\ref{1217-2}) and (\ref{1216-4}),  we obtain
\begin{eqnarray}
&&(1-|z|^2)|u(z)|\frac{|\vp^\p(z)z|}{\sqrt{1-|\vp(z)|^2}} \nonumber\\
&\lesssim& (1-|z|^2)|u(z)|\frac{|\langle \vp^\p(z) z,\vp(z)\rangle|+|\langle \vp^\p(z) z,\eta(z)\rangle|}{\sqrt{1-|\vp(z)|^2}}   \nonumber\\
&\lesssim& (1-|z|^2)|u(z)|\frac{  |\langle \vp^\p(z) z,\vp(z)\rangle|  +|\langle \vp^\p(z) z,\zeta(z)\rangle| }{{1-|\vp(z)|^2}}\label{1219-7}\\
&\lesssim& M_1 + (1-|z|^2)|u(z)| \frac{|\langle \vp^\p(z) z,\zeta(z)\rangle|}{1-|\vp(z)|^2}.\nonumber
\end{eqnarray}
By  (\ref{0102-1}), (\ref{1222-4}) and (\ref{1216-1}), we have
\begin{eqnarray}
&&(1-|z|^2)|u(z)| \frac{|\langle \vp^\p(z) z,\zeta(z)\rangle|}{1-|\vp(z)|^2} \nonumber\\
&\approx& k (1-|z|^2)|u(z)|  |\vp(z)|^{2(k-1)} |\langle \vp^\p(z) z,\zeta(z)\rangle|   \nonumber\\
&=& k (1-|z|^2)|u(z)|  |\langle \vp(z),\zeta(z)\rangle|^{k-1} |\langle \vp^\p(z) z,\zeta(z)\rangle|  \label{1219-8}\\
&\lesssim & M_1.\nonumber
\end{eqnarray}
Therefore, we have
\begin{equation}\label{1216-5}
\sup_{|\vp(z)|> \frac{1}{2}} \frac{(1-|z|^2)|u(z)|(1-|\vp(z)|^2)^\frac{1}{2}|\vp^\p(z)z|}{1-|\vp(z)|^2} \lesssim M_1.
\end{equation}
By (\ref{1216-2}, (\ref{1216-4}) and (\ref{1216-5}), we get (\ref{0101-1}).
From (\ref{1216-3}) and (\ref{0101-1}), we see that ({\it iii}) holds.

{\bf ({\it iii})$\Rightarrow$ ({\it i})}. Let $f\in H^\infty$.
From (\ref{1222-3}) and (\ref{1222-1}), we have
\begin{eqnarray}
&&\|uC_\vp f\|_{\beta} \nonumber\\
&\leq&(1-|z|^2) |\R u(z)||f(\vp(z))|  + (1-|z|^2)|u(z)|\left|(\nabla f)(\vp(z))\vp^\p(z)z\right|  \nonumber\\
&\leq&  \|u\|_{\beta}\|f\|_{\infty} +(1-|z|^2)|u(z)| \sqrt{H_{\vp(z)}(\vp^\p(z)z,\vp^\p(z)z)}\frac{\left|(\nabla f)(\vp(z))\vp^\p(z)z\right|}{ \sqrt{H_{\vp(z)}(\vp^\p(z)z,\vp^\p(z)z)}}   \nonumber\\
&\lesssim& \|u\|_{\beta}\|f\|_{\infty} +M_2 \|f\|_{\beta} \nonumber\\
&\lesssim& \|u\|_{\beta}\|f\|_{\infty} +M_2 \|f\|_{\infty}. \nonumber
\end{eqnarray}
So $uC_\vp:H^\infty\to\B$ is bounded.

{\bf ({\it iv})$\Rightarrow$ ({\it i})}. Suppose $(iv)$ holds.
For all $w\in\BB$, by Lemma 1, we have
\begin{eqnarray*}
\sqrt{\frac{2}{n+1}}\sqrt{H_{\vp(w)}(\vp^\p(w) w,\vp^\p(w) w)}
&=&\left|\phi_{\vp(w)}^\p (\vp(w))\vp^\p(w)w\right|    \\
&=&  \left(\sum_{i=1}^n \left|\R(\phi_{\vp(w),i}\comp \vp)(w)\right|^2\right)^\frac{1}{2}\\
&\lesssim &\sum_{i=1}^n \left|\R(\phi_{\vp(w),i}\comp \vp)(w)\right|.
\end{eqnarray*}
Since $|\phi_{\vp(w)}(z)|<1$, we  obtain
\begin{eqnarray}
&&(1-|w|^2)|u(w)|\sqrt{H_{\vp(w)}(\vp^\p(w) w,\vp^\p(w) w)} \nonumber\\
&\lesssim&  \sum_{i=1}^n \left (1-|w|^2)|u(w)| |\R(\phi_{\vp(w),i}\comp \vp)(w)\right|  \nonumber\\
&\lesssim& \sum_{i=1}^n\left(\|uC_\vp\phi_{\vp(w),i}\|_{\mu}+(1-|w|^2)|\R u(w)||(\phi_{\vp(w),i}\comp \vp)(w)|\right)\label{1220-1}\\
&\lesssim& M_3+\|u\|_{\mu}.\nonumber
\end{eqnarray}
Then ({\it iii}) holds. So ({\it i}) holds.

{\bf ({\it i})$\Rightarrow$ ({\it iv})}.  This implication is also obvious since $\|\phi_{\vp(w),i}\|_{\infty}\leq 1$.

{\bf ({\it v})$\Leftrightarrow$ ({\it ii})}. By $u\in\B$ and Lemma 2, we have
$$(1-|z|^2)|\R (J_u(f_{k,\xi}\comp \vp))(z) |\leq (1-|z|^2)|\R u(z)|\leq \|u\|_{\mu}$$
and
\begin{equation} \label{1223-1}
\R (uC_\vp f_{k,\xi})(z) = \R (I_u(f_{k,\xi}\comp \vp))(z)  +  \R (J_u(f_{k,\xi}\comp \vp))(z).
\end{equation}
Using triangle inequality, we can get the desired result. The proof is complete.
\end{proof}

 \noindent{\bf Remark 1.} In \cite{ZmzChh2009ams}, the equivalence of $\it (i)$ and $\it (iii)$ was proved in a different way.\msk

\section{the essential norm of $uC_\vp:H^\infty\to\B$}
To study the essential norm of $uC_\vp:H^\infty\to\B$, we need the following lemmas.\msk

 \noindent{\bf Lemma 3.} \cite[Lemma 2.10]{Tm1996d}    {\it  Suppose $T:H^\infty\rightarrow \B$ is linear and bounded. Then $T$ is  compact  if and only if whenever  $\{f_k\}_{k=1}^\infty$   is bounded in $H^\infty$ and $f_k \rightarrow 0$  uniformly  on compact subsets of $\BB$,   $\lim\limits_{k\to \infty}\|T f_{k }\|_{\B}=0$.}

\noindent{\bf Lemma 4.} {\it
Suppose $0<r,s<1$ and $f\in H(\BB)$. For all $|z|\leq s$,
$$|\nabla f(z)|\leq \frac{2n}{1-s}\max_{|z|\leq \frac{1+s}{2}}|f(z)|\,\,
\mbox{ and }\,\,
|f(z)-f(rz)|\leq \frac{2n(1-r)}{1-s}\max_{|z|\leq \frac{1+s}{2}}|f(z)|.$$
}
\begin{proof}
Since $f\in H(\BB)$, $\frac{\partial f}{\partial z_j}\in H(\BB)(j=1,2,\cdots, n)$.
For all $|z|\leq s$, let $\Gamma_{z,j,s}=\{\eta\in\DD;|\eta-z_j|=\frac{1-s}{2}\}$, and
$$z_{j,s}(\eta)=(z_1,\cdots,z_{j-1},\eta,z_{j+1},\cdots,z_n),  \eta\in\Gamma_{z,j,s}.$$
Since $|z_{j,s}(\eta)|\leq \frac{1+s}{2}$, we get
\begin{eqnarray*}
\left|\frac{\partial f}{\partial z_j}\right|
=\frac{1}{2\pi}\left|\int_{\Gamma_{z,j,s}}\frac{f(z_{j,s}(\eta))}{(\eta-z_j)^2}d\eta\right|
\leq\frac{2}{1-s}\max_{|z|\leq \frac{1+s}{2}}|f(z)|.
\end{eqnarray*}
Hence $|\nabla f(z)|\leq \frac{2n}{1-s}\max\limits_{|z|\leq \frac{1+s}{2}}|f(z)|$.  When $|z|<s$,
\begin{eqnarray*}
|f(z)-f(rz)|
&=&\left|\int_r^1 \frac{df(tz)}{dt}dt\right|
=\left|\int_r^1 \langle (\nabla f)(tz),\overline{z}\rangle dt\right| \\
&\leq& (1-r) \sup_{|z|\leq s} |\nabla f(z)|
\leq \frac{2n(1-r)}{1-s}\max_{|z|\leq \frac{1+s}{2}}|f(z)|.
\end{eqnarray*}
The proof is complete.
\end{proof}\msk

\noindent{\bf Theorem 2.} {\it Suppose $u\in H(\BB)$ and $\vp$ is a holomorphic self-map of $\BB$.  If  $uC_\vp:H^\infty\to\B$ is bounded, then 
$$   \|uC_\vp\|_{e,H^\infty\to \B}\approx Q_1+Q_2\approx Q_1+Q_3   \approx Q_1+Q_4 \approx Q_5+Q_6 .$$
Here
$$ Q_1=\limsup_{|\vp(z)|\to 1} (1-|z|^2)|\R u(z)|,\,\,Q_2= \limsup_{|\vp(z)|\to 1}(1-|z|^2)|u(z)|\sqrt{H_{\vp(z)}(\vp^\p(z) z,\vp^\p(z) z) },
$$
$$ Q_3=\limsup_{k\to\infty} \sup_{\xi\in\partial\BB}  \|u  \langle \vp,\xi\rangle^k \|_{\beta},
\,\,\,\,Q_4=\limsup_{|\vp(w)|\to 1}  \sum_{i=1}^n \|uC_\vp\phi_{\vp(w),i}\|_{\beta},
$$
$$Q_5= \limsup_{k\to\infty}\sup_{\xi\in\partial \BB}\|I_u \langle \vp,\xi\rangle^k \|_{\beta},\,\,\,\,
Q_6=\limsup_{k\to\infty}\sup_{\xi\in\partial \BB}\|J_u \langle \vp,\xi\rangle^k  \|_{\beta}.
$$
}

\begin{proof} Since $uC_\vp:H^\infty\to\B$ is bounded, by Theorem 1, we have $u\in\B$ and $\max\limits_{1\leq i\leq 6} Q_i<\infty$.

When $\sup\limits_{z\in\BB}|\vp(z)|<1$, it is easy to see that  $uC_\vp:H^\infty\to\B$ is compact by Lemmas 3 and 4.  In this case, these asymptotic relations vacuously hold. Hence we only consider the  case $\sup_{z\in\BB} |\vp(z)|=1.$

First we prove that
$$Q_1+Q_2\gtrsim \|uC_\vp\|_{e,H^\infty\to \B} .$$
Let $f_t(z)=f(tz)$, for $ f\in H(\BB)$ and $ t\in (0,1)$.  Suppose $r,s\in (\frac{1}{2},1)$.
For any $f\in H^\infty$ with  $\|f\|_{\infty}\leq 1$, we have
$$
\|(uC_\vp -uC_{r\vp})f\|_{\beta}
\leq\sup_{|\vp(z)|\leq s} G_1 +\sup_{ s<|\vp(z)|<1} G_1 +\sup_{|\vp(z)|\leq s} G_2 +\sup_{s<|\vp(z)|<1} G_2 ,
$$
where
$$G_1=(1-|z|^2)|\R u(z)||(f\comp\vp-f_r\comp\vp)(z)|,\,\,\,\,G_2=(1-|z|^2)|u(z)||\R(f\comp \vp-f_{r}\comp\vp)(z)|.$$
By Lemma 4, we have
\begin{equation}\label{1219-1}
\sup_{|\vp(z)|\leq s} G_1 \leq \frac{2n(1-r)}{1-s}\|u\|_{\beta},\,\,\sup_{s<|\vp(z)|<1} G_1 \le 2 \sup_{s<|\vp(z)|<1} (1-|z|^2)|\R u(z)|,
\end{equation}
and
\begin{eqnarray*}
\sup_{|\vp(z)|\leq s}G_2
&\leq&  \sup_{|\vp(z)|\leq s}(1-|z|^2)|u(z)||(\nabla (f-f_r))(\vp(z))|  |\vp^\p(z)z| \\
&\leq& \frac{2n}{1-s}\sup_{|\vp(z)|\leq \frac{1+s}{2}}|f(\vp(z))-f(r\vp(z)|  \sup_{|\vp(z)|\leq s} (1-|z|^2)|u(z)| |\vp^\p(z)z| \\
&\leq& \frac{8n^2(1-r)}{(1-s)^2}\sup_{|\vp(z)|\leq s} (1-|z|^2)|u(z)|  |\vp^\p(z)z| .
\end{eqnarray*}
From (\ref{1220-2}), we have
\begin{equation}\label{1219-2}
\sup_{|\vp(z)|\leq s}G_2\lesssim \frac{1-r}{(1-s)^2}M_1.
\end{equation}
By (\ref{1222-3}) and (\ref{1222-1}),  we have
\begin{eqnarray}
&&\sup_{s<|\vp(z)|<1}G_2\nonumber\\
&\leq& \sup_{s<|\vp(z)|<1}(1-|z|^2)|u(z)| \sqrt{H_{\vp(z)}(\vp^\p(z)z,\vp^\p(z)z) }\frac{\left|(\nabla f)(\vp(z))\vp^\p(z)z\right|}{ \sqrt{H_{\vp(z)}(\vp^\p(z)z,\vp^\p(z)z)} }  \nonumber\\
&&+\sup_{s<|\vp(z)|<1}(1-|z|^2)|u(z)| \sqrt{H_{\vp(z)}(\vp^\p(z)z,\vp^\p(z)z)} \frac{\left|(\nabla f_r)(\vp(z))\vp^\p(z)z\right|}{ \sqrt{H_{\vp(z)}(\vp^\p(z)z,\vp^\p(z)z)} }  \nonumber\\
&\lesssim&  (\|f\|_{\beta}+\|f_r\|_{\beta})\sup_{s<|\vp(z)|<1}(1-|z|^2)|u(z)| \sqrt{H_{\vp(z)}(\vp^\p(z)z,\vp^\p(z)z) }\nonumber
\end{eqnarray}
\begin{eqnarray}
&\lesssim& \sup_{s<|\vp(z)|<1}(1-|z|^2)|u(z)| \sqrt{H_{\vp(z)}(\vp^\p(z)z,\vp^\p(z)z) }.\label{1219-4}
\end{eqnarray}
It is obvious that
\begin{eqnarray}\label{1219-5}
\lim_{r\to 1}|(uC_{\vp}-uC_{r\vp})(0)|=0.
\end{eqnarray}
Letting $r\to 1$, by  (\ref{1219-1})-(\ref{1219-5}), we have
\begin{eqnarray*}
&&\|uC_\vp\|_{e,H^\infty\to\B}\\
&\leq& \limsup_{r\to 1} \|uC_\vp-uC_{r\vp}\|_{H^\infty\to\B}\\
&\lesssim&  \sup_{s<|\vp(z)|<1}(1-|z|^2)|u(z)| \sqrt{H_{\vp(z)}(\vp^\p(z)z,\vp^\p(z)z) }
+\sup_{s<|\vp(z)|<1} (1-|z|^2)|\R u(z)|.
\end{eqnarray*}
Here we used the fact that $uC_{r\vp}:H^\infty\to\B$ is compact.  Letting $s\to 1$, we get the desired result.

Next we prove that $$Q_1+Q_3\gtrsim Q_1+Q_2.$$
Similar to the proof of Theorem 1, we assume that there exists $k\in\N$ such that $k\geq 3 $ and $1-\frac{1}{k-1}\leq |\vp(z)| < 1-\frac{1}{k}$.
From (\ref{1219-6}), we have
\begin{eqnarray*}
\frac{|\langle  \vp^\p(z)z,\vp(z) \rangle|}{1-|\vp(z)|^2}
\lesssim \sup_{\xi\in\partial \BB}k {|\langle \vp(z),\xi \rangle^{k-1}||\langle  {\vp^\p(z)z,\xi \rangle|}}.
\end{eqnarray*}
From (\ref{1219-7}) and (\ref{1219-8}), we have
\begin{equation*}
(1-|z|^2)|u(z)|\frac{|\vp^\p(z)z|}{(1-|\vp(z)|^2)^\frac{1}{2}}
\lesssim \sup_{\xi\in\partial\BB} k(1-|z|^2)|u(z)|  |\langle \vp(z),\xi \rangle|^{k-1} |\langle \vp^\p(z) z,\xi\rangle|.
\end{equation*}
 From (\ref{1216-2}), we have
\begin{eqnarray*}
&&(1-|z|^2)|u(z)|\sqrt{H_{\vp(z)}(\vp^\p(z) z,\vp^\p(z) z)}\\
&\lesssim& \sup_{\xi\in\partial\BB} k (1-|z|^2)|u(z)|  |\langle \vp(z),\xi\rangle|^{k-1} |\langle \vp^\p(z) z,\xi\rangle|  \\
&\leq& \sup_{\xi\in\partial\BB}  \|u  \langle \varphi,\xi\rangle^k\|_{\beta}+(1-|z|^2)|(\R u)(z)|.
\end{eqnarray*}
By letting $|\vp(z)|\to 1$, we have $Q_1+Q_3\gtrsim Q_2$, i.e.,  we get $Q_1+Q_3\gtrsim Q_1+Q_2.$

Now we prove that
$$\|uC_\vp\|_{e,H^\infty\to \B} \gtrsim Q_1+Q_3.$$
 Suppose $K:H^\infty\to\B$ is compact. For any $\varepsilon>0$, there exists $\{\xi_k\}_{k=0}^\infty \subset\partial\BB$, such that
  $$ \|u\langle \varphi, \xi_k\rangle^k\|_{\beta} \geq\sup_{\xi\in\partial\BB}  \|u\langle \varphi, \xi\rangle^k\|_{\beta}-\varepsilon.$$
 Let $h_{k}(z)=\langle z,\xi_k\rangle^k.$  Thus $\|h_{k}\|_{\infty}\leq 1$ and converges to $0$ uniformly on compact subsets of $\BB$.
 Then $\lim\limits_{k\to\infty} \|Kh_k\|_{\B}=0$.  Since
\begin{eqnarray*}
\|uC_\vp -K\|_{H^\infty\to\B }
&\geq& \|(uC_\vp-K)h_{k}\|_{\B }   \\
&\geq& \|uC_\vp h_{k}\|_{\beta}+|u(0)h_{k}(\vp(0))|-\|Kh_{k}\|_{\B },
\end{eqnarray*}
we have
$$\|uC_\vp -K\|_{H^\infty\to\B}
\geq \limsup_{k\to\infty} \|uC_\vp h_{k}\|_{\beta}
\geq \limsup_{k\to\infty}\sup_{\xi\in\partial\BB}  \|u\langle \varphi, \xi\rangle^k\|_{\beta}-\varepsilon .$$
Because $K$ and $\varepsilon$ are arbitrary, we obtain  $\|uC_\vp\|_{e,H^\infty\to \B }\geq Q_3.$

Let $\{\eta_k\}_{k=1}^\infty\subset\partial\BB$ such that $Q_1=\lim\limits_{k\to \infty} (1-|\eta_k|^2)|\R u(\eta_k)|$ and $\lim\limits_{k\to\infty} |\vp(\eta_k)|=1$. Let
$$g_k(z)=\frac{2(1-|\vp(\eta_k)|^2)}{1-\langle z,\vp(\eta_k) \rangle}-\left(\frac{1-|\vp(\eta_k)|^2}{1-\langle z,\vp(\eta_k) \rangle}\right)^2.$$
Then $\{g_k\}$ is bounded in $H^\infty$ and converges to 0 uniformly on compact subsets of $\BB$. Moreover, we have
$g_k(\vp(\eta_k))=1$ and $(\R g_k)(\vp(\eta_k))=0.$ If $K:H^\infty\to\B$ is compact, we have
\begin{eqnarray*}
\|uC_\vp-K\|_{H^\infty\to\B}
&\gtrsim& (1-|\eta_k|^2)|(\R (uC_\vp g_k))(\eta_k)|-\|Kg_k\|_{\B } \\
&=&(1-|\eta_k|^2)|\R u(\eta_k)|-\|Kg_k\|_{\B }.
\end{eqnarray*}
Letting $k\to\infty$, we have $\|uC_\vp-K\|_{H^\infty\to \B}\gtrsim Q_1.$ Since $K$ is arbitrary, we get
$$\|uC_\vp\|_{e,H^\infty\to \B}\gtrsim Q_1,$$ as desired.

Next we prove that $$ \|uC_\vp\|_{e,H^\infty\to \B}\approx Q_1+Q_4.$$
  Suppose $K:H^\infty\to\B$ is compact and $1\leq i\leq n$. Let $\{w_{k}\}_{k=1}^\infty\subset\BB$
such that $\lim\limits_{k\to\infty}|\vp(w_k)|=1$  and
$$\lim\limits_{k\to\infty}\|uC_\vp \phi_{\vp(w_k),i}\|_{\mu}=\limsup_{|\vp(w)|\to 1}\|uC_\vp \phi_{\vp(w),i}\|_{\beta}.$$
By Lemma 1, $\{\phi_{\vp(w_k),i}-\vp_i(w_k) \}$ is bounded in $H^\infty$ and converges to 0 uniformly on compact subset of $\BB$. By Lemma 3,
\begin{eqnarray*}
\|uC_\vp-K\|_{H^\infty\to\B}
&\gtrsim& \limsup_{k\to\infty} \|(uC_\vp-K)(\phi_{\vp(w_k),i}-\vp_i(w_k))\|_{\beta} \\
&\geq& \limsup_{k\to\infty} \|uC_\vp\phi_{\vp(w_k),i}\|_{\beta} -\limsup_{k\to\infty} \|uC_\vp\vp_i(w_k)\|_{\beta}\\
 &&-\limsup_{k\to\infty} \|K(\phi_{\vp(w_k),i}-\vp_i(w_k))\|_{\beta}  \\
 &\geq& \limsup_{|\vp(w)|\to 1}\|uC_\vp \phi_{\vp(w),i}\|_{\beta}-Q_1.
\end{eqnarray*}
Since $K$ is arbitrary, we have $\|uC_\vp\|_{e,H^\infty\to \B}\gtrsim Q_4 -Q_1$. Since $\|uC_\vp\|_{e,H^\infty\to \B}\gtrsim Q_1$, we get $$\|uC_\vp\|_{e,H^\infty\to \B}\gtrsim Q_1+Q_4.$$
 From (\ref{1220-1}), we have  $Q_1+Q_4\gtrsim Q_2$. So $$Q_1+Q_4\gtrsim Q_1+Q_2\gtrsim \|uC_\vp\|_{e,H^\infty\to \B},$$ as desired.

Finally, we prove that $$ Q_5+Q_6 \approx Q_3+Q_1.$$

From (\ref{1223-1}), we have $Q_3\leq Q_5+Q_6,\,\,\mbox{ and }\,\,Q_5\leq Q_6+Q_3.$
By Lemma 2,  $Q_6\leq Q_1$. So $$Q_5+Q_6\lesssim Q_1+Q_3.$$

Suppose $\{z_k\}_{k=1}^\infty\subset\BB$ such that $\lim\limits_{k\to\infty} (1-|z_k|)|(\R u)(z_k)|=Q_1$ and $\lim\limits_{k\to\infty}|\varphi(z_k)|=1$.  From the fact that 
\begin{eqnarray*}
(1-|z_k|)|(\R u)(z_k)|
&\approx &(1-|z_k|)|(\R u)(z_k) |\vp(z_k)|^{2k}\\
&\leq &\|J_u f_{k,\vp(z_k)} \|_{\beta}
\leq \sup_{\xi\in\partial\BB}\|J_u  \langle\varphi, \xi\rangle^k\|_{\beta},
\end{eqnarray*}
we have $Q_1\lesssim Q_6$.  So $ Q_1+Q_3\lesssim Q_5+Q_6.$ The proof is complete.
\end{proof}

From Theorem 2, we immediately get the following corollary.\msk

\noindent{\bf Corollary 1.} {\it Suppose $u\in H(\BB)$ and $\vp$ is a holomorphic self-map of $\BB$.  If  $uC_\vp:H^\infty\to\B$ is bounded, then the following statements are equivalent.

\begin{enumerate}[(i)]
  \item $uC_\vp:H^\infty\to\B$ is compact.
  \item $$\limsup_{|\vp(z)|\to 1} (1-|z|^2)|\R u(z)|=0 ~~~\mbox{and}~~~ \limsup_{k\to\infty} \sup_{\xi\in\partial\BB}  \|u  \langle \vp,\xi\rangle^k \|_{\beta}=0.$$
  \item $$\limsup_{|\vp(z)|\to 1} (1-|z|^2)|\R u(z)|=0 ~~~\mbox{and}~~~  \limsup_{|\vp(z)|\to 1}(1-|z|^2)|u(z)|\sqrt{H_{\vp(z)}(\vp^\p(z) z,\vp^\p(z) z) }=0.$$
  \item $$\limsup_{|\vp(z)|\to 1} (1-|z|^2)|\R u(z)|=0 ~~~\mbox{and}~~~  \limsup_{|\vp(w)|\to 1}  \sum_{i=1}^n \|uC_\vp\phi_{\vp(w),i}\|_{\beta}=0.$$
  \item $$ \limsup_{k\to\infty}\sup_{\xi\in\partial \BB}\|I_u \langle \vp,\xi\rangle^k \|_{\beta}=0 ~~\mbox{and}~~~
 \limsup_{k\to\infty}\sup_{\xi\in\partial \BB}\|J_u \langle \vp,\xi\rangle^k  \|_{\mu}=0.
$$
\end{enumerate}
}

\noindent{\bf Remark 2.}   Suppose $u,v\in H(\BB)$ and $\vp,\psi$ are  holomorphic self-maps of $\BB$.  Based on the work of \cite{SycLsxZxl2007}, we  conjecture that the following statements hold:

{\it (a)}   $uC_\vp-vC_\psi:H^\infty\to\B$ is bounded if and only if
  $$\sup_{k\in \N\cup\{0\}} \sup_{\xi\in\partial \BB} \|(uC_\vp-vC_\psi)\langle z,\xi\rangle^k\|_\beta <\infty.$$

{\it (b)}  Assume that $uC_\vp :H^\infty\to\B$ and $ vC_\psi:H^\infty\to\B$  are bounded. Then
\begin{eqnarray*}
\|uC_\varphi-vC_{\psi}\|_{e,H^\infty\to\B }\thickapprox \limsup_{k\to\infty} \sup_{\xi\in\partial \BB}\|(uC_\vp-vC_\psi)\langle z,\xi\rangle^k\|_\beta .\nonumber
\end{eqnarray*}
We are not able, at the moment, to prove this conjecture. Hence, we leave the problem to the readers
interested in this research area.


\begin{thebibliography}{aa}
\bibitem{Cf2013cejm} F. Colonna, New criteria for boundedness and compactness of weighted composition operators mapping into the Bloch space, {\it Cent. Eur. J. Math.} \textbf{11} (2013), 55--73.

\bibitem{CccMbd1995} C.~C.~Cowen and  B.~D.~MacCluer, {\it Composition Operators on Spaces of Analytic Functions},  CRC Press, Boca Raton, FL,
1995.

\bibitem{Djn2012jmaa} J. Dai, Compact composition operators on the Bloch space of the unit ball, {\it J. Math. Anal. Appl.} {\bf 386} (2012), 294--299.

\bibitem{HqhLsxWh2017cvee} Q. Hu, S. Li and H. Wulan, New essential norm estimates of weighted composition operators from $H^\infty$ into the Bloch space, {\it Complex Var. Elliptic Equ.} {\bf 62} (2017), 600--615.

\bibitem{li} S.~Li, Riemann-Stieltjes operators from $F(p,q,s)$ to Bloch space on the unit ball, {\it J.
Inequal. Appl.} Vol. 2006 (2006), Article ID 27874, 14 pages.

\bibitem{LsxSs2007aaa} S.~Li and S.~Stevi\'c, Weighted composition operators from $H^\infty$ to the Bloch space on the polydisc, {\it Abstr. Appl. Anal.} Vol. 2007 (2007), Article ID 48478, 12 pages.

\bibitem{LsxSs2008tjm} S.~Li and S.~Stevi\'c, Weighted composition operators between $H^\infty$ and $\alpha$-Bloch spaces in the
unit ball,  {\it Taiwanese J. Math.} {\bf 12} (2008), 1625--1639.

 \bibitem{LxsLsx2017ieot} X. Liu and S.~Li, Norm and essential norm of a weighted composition operator on the Bloch space, {\it Integr. Equ. Oper. Theory} {\bf 87} (2017), 309--325.

\bibitem{MkMa1995tams}  K.~Madigan and A.~Matheson, Compact composition operators on the Bloch space, {\it Trans. Amer. Math. Soc.} {\bf 347} (1995), 2679--2687.

\bibitem{Os2001tjm} S. Ohno, Weighted composition operators between $H^\infty$ and Bloch space. {\it Taiwanese J. Math.} {\bf 5} (2001), 555--563.

\bibitem{SycLsxZxl2007} Y. Shi, S. Li and X. Zhu,  Differences of weighted composition operators from $H^\infty $ to the Bloch space,  {\it  arXiv:1712.03402} (2017), 18 pages.

\bibitem{SjhLl2000ams} J.~Shi and L.~Luo, Composition operators on the Bloch space, {\it Acta Math. Sin.} {\bf 16} (2000), 85--98.

\bibitem{Tr1980blms} R. Timoney,  Bloch function in several complex variables, I,  {\it Bull. London Math. Soc.} {\bf 12} (1980), 241--267.

\bibitem{Tm1996d} M. Tjani, {\it Compact composition operators on some M\"obius invariant Banach spaces}, PhD dissertation, Michigan State University, 1996.

\bibitem{WlZdcZkh2009pams} H. Wulan, D. Zheng and K. Zhu, Compact composition operators on BMOA and the Bloch space, {\it Proc. Amer. Math. Soc.}  {\bf 137} (2009), 3861--3868.

\bibitem{ZmzChh2009ams} M. Zhang and H. Chen, Weighted composition operators of $H^\infty$ into $\alpha$-Bloch spaces on the unit ball, {\it Acta Math. Sin.} {\bf 25} (2009), 265--278.

\bibitem{Zkh2005} K. Zhu, {\it Spaces of Holomorphic Functions  in the Unit Ball}, Springer, New York, 2005.

\end{thebibliography}
\end{document}